\documentclass[12pt]{amsart}
\usepackage{amsmath,amscd, xcolor, setspace, amsfonts, bbm, amsthm, tikz, tikz-cd, enumerate, graphicx}
\usetikzlibrary{calc,decorations.markings}
\usepackage[margin=1in]{geometry}
\usepackage[shortalphabetic]{amsrefs}
\allowdisplaybreaks
\usepackage{setspace}
\setdisplayskipstretch{2.5}
\newtheorem{thm}{Theorem}[section]
\newtheorem{theorem}[thm]{Theorem}

\newtheorem{lemma}[thm]{Lemma}
\newtheorem*{lemma*}{Lemma}

\newtheorem{cor}[thm]{Corollary}

\newtheorem*{remark*}{Remark}
\newtheorem*{grhthm*}{Theorem \ref{GRHthm}}
\newtheorem*{chithm*}{Theorem \ref{typicalchi}}
\numberwithin{equation}{section}

\newcommand{\Q}{\mathbb{Q}}
\newcommand{\Z}{\mathbb{Z}}

\newcommand{\C}{\mathbb{C}}

\newcommand{\res}{\underset{\small s=1}{\mathrm{Res}}~}
\title{Comparing the Density of $D_4$ and $S_4$ Quartic Extensions of Number Fields}
\author{Matthew Friedrichsen and Daniel Keliher}

\begin{document}
\maketitle

\begin{abstract}
   When ordered by discriminant, it is known that about 83\% of quartic fields over $\mathbb{Q}$ have associated Galois group $S_4$, while the remaining 17\% have Galois group $D_4$. We study these proportions over a general number field $F$. We find that asymptotically 100\% of quadratic number fields have more $D_4$ extensions than $S_4$ and that the ratio between the number of $D_4$ and $S_4$ quartic extensions is biased arbitrarily in favor of $D_4$ extensions. Under GRH, we give a lower bound that holds for general number fields.
\end{abstract}

\section{Introduction and Statement of Main Results}

From Hilbert's Irreducibility Theorem it follows that the Galois group of the splitting field of a ``random'' degree $n$ polynomial over $\Q$ will be $S_n$ 100\% of the time \cite{Serre}. If we instead pick a random degree $n$ extension of $\Q$ ordered by discriminant, one might expect the same behavior. Indeed, this is clearly true for $n = 2$, and for $n=3$ this is true by work of Davenport and Heilbronn \cite{DavHeil}. For $n=4$, work of Bhargava \cite{Bharg} and Cohen, Diaz y Diaz, and Olivier \cite{CDO} shows that only about 83\% of quartic extensions of $\Q$ have Galois group $S_4$, with the remaining 17\% having Galois group $D_4$ and 0\% having Galois groups $C_4, V_4,$ or $A_4.$ In our work, we investigate this disparity for quartic extensions of an arbitrary number field $F$. In particular, we ask what proportion of quartic extensions of $F$ are $S_4$ and what proportion are $D_4.$

Our first result shows that, when $F$ is quadratic, there are typically many more $D_4$ than $S_4$ extensions. To make this precise, let $N_n^F(X;G) = \#\{K/F : |D_{K/F}| < X, [K:F] = n, \text{Gal}(\widetilde{K}/F) = G\}.$

\begin{theorem} \label{maintheorem}
For $\epsilon> 0$, asymptotically 100\% of quadratic number fields $F$ ordered by discriminant have
$$\lim_{X \to \infty} \frac{N_4^F(X;D_4)}{N_4^F(X;S_4)} \gg (\log |D_F|)^{\log 2 - \epsilon}.$$
\end{theorem}

In particular, we have:

\begin{cor}\label{mored_4}
$100\%$ of quadratic number fields $F$ have arbitrarily many more $D_4$ 	quartic extensions than $S_4$ quartic extensions. 
\end{cor}

In practice, one can find quadratic number fields with small discriminant where $D_4$ quartic extensions vastly outnumber $S_4$ quartic extensions. For example, more than 90\% of quartic extensions of $\Q(\sqrt{-210})$ and more than 99\% of quartic extensions of $\Q(\sqrt{-510510})$ are $D_4$ quartic extensions. We give a general lower bound on the ratio for quadratic number fields in Theorem \ref{nonstat}.

For general number fields we prove the following conditional statement. Let $\mathrm{Cl}_F$ be the the ideal class group of a number field $F$ and let $\#\mathrm{Cl}_F[2]$ be the number of elements of $\mathrm{Cl}_F$ with order dividing 2. Then we have the following:

\begin{theorem}\label{GRHthm}
Assume GRH and let $F$ be a degree $d$ number field over $\Q$. Then, 
$$\lim_{X \rightarrow \infty}\frac{N^F_4(X; D_4)}{N^F_4(X;S_4)} \gg_d \frac{\#\mathrm{Cl}_F[2]-1}{(\log \log |D_F|)^d}.$$
\end{theorem}
The assumption of GRH is used for a lower bound on the residue of Dedekind zeta functions. In the course of proving Theorem \ref{maintheorem}, we prove that a weaker, but still sufficient bound holds for a positive proportion of quadratic Dirichlet $L$-functions in a restricted family. In particular, we show:

\begin{theorem}\label{typicalchi}
For 100\% of fundamental discriminants $D$ and for $\epsilon,\delta > 0,$ a proportion $1-\delta$ of quadratic characters $\chi \pmod{|D|}$ have
$$L(1,\chi) \ge \exp \left(-c(\log\log|D|)^{1-\frac{\log 2}{2} + \epsilon} \right),$$
where $c$ depends on $\delta$.
\end{theorem}

Granville and Soundararajan in \cite{GranvilleSound} study the distribution of $L(1, \eta_D)$, where $\eta_D$ is the primitive real character with modulus $|D|$, as $D$ ranges over fundamental discriminants with $|D| \le x$. Since we need to restrict our attention to the typical behavior of $L(1, \chi)$ for the much smaller family of quadratic characters of a fixed modulus $|D|$, their results do not port over directly to this setting.

In the next section we will show how to use field counting results of Bhargava, Shankar and Wang and of Cohen, Diaz y Diaz and Olivier \cite{BSW,CDO} to prove Theorem \ref{GRHthm}. In the following section, we consider the family of quadratic Dirichlet $L$-functions and prove Theorem \ref{typicalchi}. In Section 4, we complete the proof of Theorem \ref{maintheorem}, and in Section 5, we provide some examples.

\section{Field Counting and Proof Strategy}

For an extension of number fields $L/F$, let $(r_1, r_2)$ denote the signature of $F$, $D_F$ the absolute discriminant of $F$, and $D_{L/F}$ the norm of the relative discriminant of $L/F$. Note that $D_L=D_{L/F} D_F^{[L:F]}$. As above, $\mathrm{Cl}_F$ denotes the ideal class group of $F$ and $\mathrm{Cl}_F[2]$ the elements of $\mathrm{Cl}_F$ with order dividing 2. 

Bhargava, Shankar, and Wang \cite{BSW} give asymptotic formulas for $N^F_n(X; S_n)$ when $n=2,3,4$ or $5$. In the $n=4$ case they prove:
\begin{theorem}[Bhargava, Shankar, Wang]
If $F$ is a number field with $r_1$ real embeddings and $r_2$ complex embeddings, then
\begin{equation}\label{S4formula}
N_4^F(X;S_4) \sim X \frac{1}{2} ~ \res~\zeta_F(s) \left( \frac{10}{4!}\right)^{r_1}\left(\frac{1}{4!}\right)^{r_2} \prod_{\frak{p}}\left(1+\frac{1}{N\frak{p}^2} - \frac{1}{N\frak{p}^3} - \frac{1}{N\frak{p}^4}\right),
\end{equation}
where the product runs over prime ideals of $F$.
\end{theorem}
It follows that $\lim_{X\rightarrow \infty} \frac{1}{X}N_4^F(X; S_4) \asymp \res\zeta_F(s)$. Thus, it is bounds on the residue of $\zeta_F(s)$ that we'll need to control this term. For $D_4$, we recall work of Cohen, Diaz y Diaz, and Olivier \cite{CDO} that gives an asymptotic formula for $N_4^F(X;D_4)$. 
\begin{theorem}[Cohen, Diaz y Diaz, Olivier]
If $F$ is a number field with $r_2$ complex embeddings, then
\begin{equation}\label{D4formula} 
N^F_4(X; D_4) \sim X\sum_{ [L:F]=2} \frac{1}{2^{r_2+1}D_{L/F}^2 \zeta_L(2)} \res~\zeta_L(s),
\end{equation}
where the sum runs over quadratic extensions of $F$.
\end{theorem}
From these we obtain upper bounds on $N^F_4(X; S_4)$ and lower bounds on $N^F_4(X; D_4)$ so as to bound their ratio from below.   

Restricting the summation in (\ref{D4formula}) to be over only those quadratic extensions $L$ of $F$ that are unramified, i.e. $L$ such that $D_{L/F}=1$, yields, 
\begin{equation}\label{D4trivbound}
\lim_{X \rightarrow \infty}\frac{N^F_4(X; D_4)}{X} \gg_d \sum_{\substack{ [L:F]=2 \\ D_{L/F}=1}} \underset{s=1}{\text{Res}}~\zeta_L(s).
\end{equation}
If $K$ is a number field of degree $d$ over $\Q$ of discriminant $D_K$, then under GRH we have (see e.g. \cite{CK}),
\begin{equation}\label{GRHbound}
\frac{1}{\log\log |D_K|} \ll \underset{s=1}{\text{Res}}~\zeta_K(s) \ll (\log \log |D_K|)^{d-1}.
\end{equation}

Applying this bound to (\ref{S4formula}) gives an asymptotic upper bound for $N_4^F(X;S_4)$ and likewise applying it to $(\ref{D4trivbound})$ gives a lower bound for $N_4^F(X;D_4)$. In particular, conditional on GRH,
\begin{equation} \label{d4s4bounds}
\lim_{X \rightarrow \infty}\frac{N^F_4(X; S_4)}{X} \ll (\log \log |D_F|)^{d-1} ~\text{and}~ \lim_{X \rightarrow \infty}\frac{N^F_4(X; D_4)}{X} \gg \sum_{\substack{ [L:F]=2 \\ D_{L/F}=1}} \frac{1}{\log \log |D_L|}
\end{equation}
For a number field $F$, it follows from class field theory that there are $\#\mathrm{Cl}_F[2]-1$ quadratic extensions $L/F$ such that $D_{L/F}=1$. Using this fact and the estimates (\ref{d4s4bounds}) we bound the ratio $N^F_4(X; D_4)/N^F_4(X;S_4)$. We immediately obtain Theorem \ref{GRHthm}. Note that for number fields $F$ with odd class number, the lower bound given by the theorem is $0$. However, you could obtain a similar lower bound by instead summing over quadratic extensions $L/F$ with $D_{L/F}$ up to some bound.

If we specialize Theorem \ref{GRHthm} to quadratic number fields $F$ and note that $\#\mathrm{Cl}_F[2] = 2^{\omega(D_F)-m}$, where $m=1$ or $2$ and $\omega(n)$ is the number of distinct prime divisors of $n$, we obtain the following:

\begin{cor}\label{quadcor}
Assume GRH and let $F$ be any quadratic number field. Then,
$$\lim_{X \rightarrow \infty}\frac{N^F_4(X; D_4)}{N^F_4(X;S_4)} \gg \frac{2^{\omega(D_F)}}{(\log \log |D_F|)^2}.$$
\end{cor}

Even when $\#\mathrm{Cl}_F$ is odd, the lower bound from the corollary still holds because the sum on the right hand side of (\ref{d4s4bounds}) can be expanded to include extensions $L/F$ that ramify only at primes dividing $2$.

The rest of the paper is essentially concerned with removing the GRH assumption from Corollary \ref{quadcor}. To do this we will restrict our attention to the case where $F$ is a quadratic number field, and prove bounds analogous to (\ref{GRHbound}) unconditionally for a positive proportion of the necessary $L$-functions. Theorem \ref{maintheorem} will follow from such bounds.

\section{Typical behavior of $L(1,\chi)$}\label{L-functions}

We begin with a lemma that isolates the $L$-functions of interest.

\begin{lemma}\label{chifacts}
Let $F$ be a quadratic number field and $L/F$ be an unramified quadratic extension. Then there are nonprincipal quadratic Dirichlet characters $\chi_1$ and $\chi_2$ such that $\chi_1\chi_2 = \chi_F$ and for which
$$\zeta_L(s) = \zeta(s)L(s,\chi_F)L(s,\chi_1)L(s,\chi_2),$$
where $\chi_F = \left(\frac{D_F}{\cdot}\right).$
\end{lemma} 

When we put this in the context of the ratio $N_4^F(X;D_4)/N_4^F(X; S_4)$, we see that each such unramified quadratic extension $L$ of $F$ has $\res \zeta_F(s)$ dividing $\res \zeta_L(s),$ leaving behind 
\begin{equation}\label{residueratio}
\frac{\res \zeta_L(s)}{\res \zeta_F(s)} = L(1,\chi_1)L(1,\chi_2).
\end{equation}

Therefore, in order to remove the GRH assumption and prove Theorem \ref{maintheorem}, it suffices to prove Theorem \ref{typicalchi}.

Rather than bound $L(1,\chi)$ directly, we'll consider $\log L(1,\chi)$ instead. We will find its second moment and then use a discrete analogue of Chebyshev's inequality to obtain the desired result.

\subsection{The Second Moment of $\log L(1,\chi)$}

Recall from \cite[Chapter 20]{Dav} that at most one nonprincipal real character $\chi \pmod{D}$ exists such that $L(s,\chi)$ has a real zero $\beta$ with $\beta > 1 - c/\log D$ for some absolute constant $c$. If such a character exists, we say it has an exceptional zero. We will use this in the following key lemma.
\begin{lemma}\label{secmoment}
For a fixed modulus $D$, let $\chi_0$ denote the principal character and define\\ $V=\{\chi \pmod D \mid \chi\neq \chi_0, \chi^2=\chi_0,~ \chi \text{ does not have an exceptional zero} \}.$ Then
$$\frac{1}{\#V} \sum_{\chi \in V} (\log L(1,\chi))^2 = O\left((\log\log D)(\log \omega_Y(D) + \log\log Y) + \frac{(\log\log D)^2}{2^{\omega_Y(D)}} \right),$$ for any $Y < D$ and where $\omega_Y(D) = \#\{p \le Y : p \mid D \}.$
\end{lemma}
To prove the lemma, we need the following result of Brun and Titchmarsh \cite[Theorem 2]{BT}.
\begin{theorem}[Brun-Titchmarsh] \label{B-T}
Let $a$ and $q$ be coprime integers, $\pi(z; q,a)$ be the number of primes less that $z$ that are congruent to $a \pmod q$, and let  $x\geq 0$ and $y > q$ be real numbers. Then 
$$\pi(x+y;q,a)-\pi(x;q,a) \leq \frac{2y}{\varphi(q)\log y/q}.$$
\end{theorem}

\begin{proof}[Proof of Lemma \ref{secmoment}] Recall that for a character $\chi \pmod{D},$
$$\log L(1,\chi) = \sum_{p} \frac{\chi(p)}{p} + \sum_{p}\sum_{k=2}^{\infty} \frac{\chi(p)^k}{kp^k}.$$
The double sum on the right is absolutely convergent and bounded above and below by absolute constants. So we need only focus on the left-hand sum. Because the left-hand sum is not absolutely convergent, we will split the sum at some threshold $T$ into
$$\sum_{p < T} \frac{\chi(p)}{p} + \sum_{p \ge T}\frac{\chi(p)}{p}.$$
From \cite[Chapter 20]{Dav}, if $\chi$ does not have an exceptional zero, a choice of $T = D^{\log D/(4c_1^2)}$ for some absolute constant $c_1$, yields
$$\sum_{p \ge T} \frac{\chi(p)}{p} = O\left(\frac{1}{\exp(c_2\sqrt{\log T})\log T} \right),$$
where $c_2$ is an absolute constant. This makes
\begin{equation}\label{bigOlog}
\log L(1,\chi) = \sum_{p < T} \frac{\chi(p)}{p} + O(1).
\end{equation}

If we assume the set of all Dirichlet characters modulo $D$ does not contain an exceptional zero and we use (\ref{bigOlog}), then
\begin{align*}
	\frac{1}{\#V}\sum_{\chi \in V} (\log L(1,\chi))^2 & = \frac{1}{\#V}\sum_{\chi \in V} \left( \sum_{p} \frac{\chi(p)}{p} + \sum_{p}\sum_{k=2}^{\infty} \frac{\chi(p)^k}{kp^k}\right)^2 \\
	&= \frac{1}{\#V}\sum_{\chi \in V} \left(\sum_{p < T} \frac{\chi(p)}{p} + O(1) \right)^2 \\
	& = \frac{1}{\#V}\sum_{\chi \in V}\left(\sum_{p,q < T} \frac{\chi(pq)}{pq} + O(1)\sum_{p < T}\frac{\chi(p)}{p} + O(1)\right).
\end{align*}

The cross term we can bound trivially as $O(\log \log T)$, so the object of our focus will be the leading term as we sum over the characters in our family. For the leading term, we include the principal character in the outer sum, apply orthogonality and then subtract off the contribution from the principal character. After summing over $\chi$, this yields a main term of

\begin{equation}\label{ortho}
\sum_{p<T} \frac{1}{p}\sum_{\substack{{q< T}\\{pq \,\equiv\, \square\,(D)}}} \frac{1}{q} - \frac{1}{\#V + 1}\sum_{p,q<T}\frac{1}{pq}.
\end{equation}

Given that $\#V + 1 = 2^{\omega(D)},$ the right-hand term of (\ref{ortho}) is on the order of $(\log \log T)^2\big/2^{\omega(D)}.$ Since our goal is to show the second moment is small, we need to calculate an upper bound for the left-hand term of (\ref{ortho}). For the analysis, we will break up the left-hand sum of (\ref{ortho}) according to an auxiliary parameter $E<D$,
\begin{equation}\label{ortho2}
\sum_{p<T} \frac{1}{p}\sum_{\substack{{q< T}\\{pq \,\equiv\, \square\,(D)}}} \frac{1}{q} = \sum_{p<T} \frac{1}{p} \Bigg(\sum_{\substack{{q< 2E}\\{pq \,\equiv\, \square\,(D)}}} \frac{1}{q} + \sum_{\substack{{2E \le q< 2D}\\{pq \,\equiv\, \square\,(D)}}} \frac{1}{q} + \sum_{\substack{{2D \le q< T}\\{pq \,\equiv\, \square\,(D)}}} \frac{1}{q}\Bigg).
\end{equation}

We will a make a convenient choice of $E$ later. For now, we will establish an upper bound on the rightmost sum over primes $2D \le q < T$ in (\ref{ortho2}). Note that given a fixed prime $p$, there are exactly $\dfrac{\phi(D)}{2^{\omega(D)}}$ congruence classes $a \in (\Z/D\Z)^{\times}$ such that $pa$ is a square modulo $D$. Hence,

\begin{equation}\label{midqs}
\sum_{\substack{2D \le q < T\\pq \equiv \square (D)}} \frac{1}{q} = \sum_{\substack{a \in (\Z/D\Z)^{\times}\\pa \equiv \square (D)}} \sum_{\substack{2D \le q < T\\ q \equiv a (D)}} \frac{1}{q}.
\end{equation}

Using partial summation, the inner sum of (\ref{midqs}) is
\begin{align*}
  \sum_{\substack{2D \le q < T\\ q \equiv a (D)}} \frac{1}{q}  & =\int_{2D}^T \frac{1}{t}\, d(\pi(t;a,D)) \\
    & = \frac{\pi(t;a,D)}{t}\Big|_{2D}^T + \int_{2D}^T \frac{\pi(t;a,D)}{t^2}\,dt.
\end{align*}
We apply Theorem \ref{B-T} to bound from above $\pi(t; a, D)$ when $t> D$.
\begin{align*}
    \frac{\pi(t;a,D)}{t}\Big|_{2D}^T + \int_{2D}^T \frac{\pi(t;a,D)}{t^2}\,dt & \le \frac{2}{\phi(D)}\left(\frac{1}{\log(T/D)} - \frac{2}{\log 2}\right) + \frac{2}{\phi(D)}\int_{2D}^T \frac{1}{t \log (t/D)}\,dt \\
 	& = \frac{2}{\phi(D)}\left(\log \log (T/D) - \log \log 2 +  \frac{1}{\log(T/D)} - \frac{2}{\log 2}\right).
\end{align*}

Note that the bound  we obtained above by using Theorem \ref{B-T} doesn't depend on $a$, so we may use the same bound for every relevant congruence class. This means the double sum (\ref{midqs}) can be bounded by
$$\sum_{\substack{a \in (\Z/D\Z)^{\times}\\pa \equiv \square (D)}} \sum_{\substack{2D \le q < T\\ q \equiv a (D)}} \frac{1}{q} \le \frac{1}{2^{\omega(D)-1}}\left(\log \log (T/D) - \log \log 2 +  \frac{1}{\log(T/D)} - \frac{2}{\log 2}\right).$$
Overall, this term is of order $O\left(\dfrac{\log\log D}{2^{\omega(D)}}\right).$

Because Theorem \ref{B-T} does not apply when $ t \le D$, we need to handle the range $2E \le q < 2D$ differently. Note that if $E \mid D$, then
$$\sum_{\substack{R \le q < S\\q \equiv a (D)}} \frac{1}{q} \le \sum_{\substack{R \le q < S\\q \equiv a (E)}}\frac{1}{q},$$
for any choice of $1 \le R < S.$

We will take $E$ to be a sufficiently small divisor of $D$.  Take $Y \le D$, let $\omega_Y(D) = \#\{p \mid D : p \le Y\}$, and define
$$E = \prod_{\substack{p \le Y\\p \mid D}} p.$$

We can now use Theorem \ref{B-T} with $E$ as the modulus to obtain
$$ \sum_{\substack{{2E \le q< 2D}\\{pq \,\equiv\, \square\,(D)}}} \frac{1}{q} \le \sum_{\substack{{2E \le q< 2D}\\{pq \,\equiv\, \square\,(E)}}} \frac{1}{q} \le \frac{1}{2^{\omega_Y(D)-1}}\left(\log \log (2D/E) - \log \log 2 +  \frac{2}{\log(2D/E)} - \frac{2}{\log 2} \right)$$
where $\omega_Y(D)$ is the number of distinct prime divisors of $D$ that are less than or equal to $Y$. Overall this is of order $O\left(\dfrac{\log\log D}{2^{\omega_Y(D)}} \right).$

Finally, for the range $q < 2E$, we will use the trivial bound on the sum of reciprocal primes,

$$\sum_{\substack{{q< 2E}\\{pq \,\equiv\, \square\,(D)}}} \frac{1}{q} \le \sum_{q < 2E} \frac{1}{q} = \log \log 2E + O(1).$$
Given that $E$ is the product of primes that are at most $Y$, then
\begin{align*}
\log \log 2E & \le \log\log (2Y^{\omega_Y(D)}) \\
& = \log (\log 2 + \omega_Y(D)\log Y) \\
& = O(\log \omega_Y(D) + \log \log Y).
\end{align*}

Now, because

$$\sum_{p < T} \frac{1}{p} = \log \log T + O(1) = O(\log \log D)$$ and $\dfrac{(\log\log D)^2}{2^{\omega_Y(D)}} \gg \dfrac{(\log\log D)^2}{2^{\omega(D)}},$ then our whole sum is

\begin{equation}\label{orderofsecond}
O\left((\log\log D)(\log \omega_Y(D) + \log\log Y) + \frac{(\log\log D)^2}{2^{\omega_Y(D)}}\right).
\end{equation}

This completes the proof of the lemma if the family of characters does not contain an exceptional zero. If the complete family of characters modulo $D$ does admit an exceptional zero, then we also need to subtract off the contribution from the exceptional character $\chi'$. Then (\ref{ortho}) becomes
$$\sum_{p<T} \frac{1}{p}\sum_{\substack{{q< T}\\{pq \,\equiv\, \square\,(D)}}} \frac{1}{q} - \frac{1}{\#V + 2}\sum_{p,q<T}\frac{1}{pq} - \frac{1}{\#V + 2}\sum_{p,q < T} \frac{\chi'(pq)}{pq}.$$
The term coming from the exceptional character is $O\left(\dfrac{(\log \log D)^2}{2^{\omega(D)}}\right)$ and does not change the rest of the analysis.
\end{proof}

Now that we know the order of the second moment for the family of quadratic Dirichlet $L$-functions, we want to use this to give a workable lower bound on $L(1,\chi)$ for most quadratic $\chi$ in $V$. This result is the following corollary.

\begin{cor}\label{generallowerL}
For a modulus $D$, any $Y\le D$, and a choice of $k \ge 1$ a proportion at least $1 - 1/k^2$ of non-exceptional quadratic characters $\chi$ modulo $D$ have
$$L(1,\chi) \ge (\log D)^{-kc\sqrt{\frac{\log \omega_Y(D) + \log \log Y}{\log \log D} + \frac{1}{2^{\omega_Y(D)}}}}, $$
where $c$ is an absolute constant.
\end{cor}

\begin{proof}
Let $k \ge 1$ and $\sigma^2$ be any value such that
$$\frac{1}{\#V} \sum_{\chi \in V} (\log L(1,\chi))^2 \le \sigma^2.$$ 
We bound $\#\{\chi \in V : |\log L(1,\chi)| \ge k\sigma\}$.
\begin{equation*}
\#\{\chi \in V : |\log L(1,\chi)| \ge k\sigma\} = \sum_{\substack{\chi \in V\\|\log L(1,\chi)| \ge k\sigma}}1 \le \sum_{\chi \in V} \frac{(\log L(1,\chi))^2}{k^2\sigma^2} \le \frac{\#V}{k^2}.
\end{equation*}
The result follows from this argument and Lemma \ref{secmoment}.
\end{proof}

\subsection{Typical Behavior of $\omega_Y$}

From (\ref{orderofsecond}), the order of the second moment depends on the size of $\omega_Y(D).$ We will apply Chebyshev's inequality to give workable bounds on $\omega_Y(D)$. Because we will make frequent use of Chebyshev's inequality, we state it here for convenience.

\begin{theorem}[Chebyshev's Inequality]\label{Cheby}
Let $V$ be a finite set with cardinality $N$ and let $f\colon V \rightarrow \C$ be function on $V$. If $f$ has mean and variance respectively given by
$$\mu = \frac{1}{N}\sum_{v \in V} f(v) \text{~~ and ~~} \sigma^2 = \frac{1}{N} \sum_{v \in V} (f(v) -\mu)^2,$$
then, for any $k$,
$$\#\{v \in V \mid |f(v)- \mu| \geq k\sigma\} \leq \frac{N}{k^2}.$$
\end{theorem}

To use this, we will need to understand the mean and variance of $\omega_Y(D)$ for most $D$.

 Let $N(X)$ be the number of fundamental discriminants less than $X$. It is known that $N(X) \sim c_2X$ for some constant $c_2$.

\begin{lemma}[Mean of $\omega_Y(n)$]\label{omegamean}
Let $Y \le X^{1 - \delta}$ for some $\delta > 0$. For fundamental discriminants $D$ such that $|D| \le X,$
$$\frac{1}{N(X)} \sideset{}{^\flat}\sum_{|D| \le X} \omega_Y(D) = \sum_{p \le Y} \frac{1}{p+1} + O\left(\frac{Y^{1/2}}{X^{1/2 - \epsilon}}\right).$$
\end{lemma} 
\begin{proof}
We have
\begin{align*}
	\frac{1}{N(X)}\sideset{}{^\flat}\sum_{D \le X} \omega_Y(D) & =  \frac{1}{N(X)}\sideset{}{^\flat}\sum_{|D| \le X} \sum_{\substack{p \mid D \\ p \le Y}} 1 \\
	& =\frac{1}{N(X)} \sum_{p \le Y}  \sideset{}{^\flat}\sum_{\substack{p \mid D \\ |D| \le X}} 1 \\
	& = \frac{1}{c_2X} \sum_{p \le Y} \left(c_2\frac{X}{p+1} + O\left(X^{1/2+\epsilon}p^{-1/2}\right)\right) \\
	& = \sum_{p \le Y} \left(\frac{1}{p+1} + O(X^{-1/2 + \epsilon}p^{-1/2})\right) \\
	& = \sum_{p \le Y} \frac{1}{p+1} + O\left(\frac{Y^{1/2}}{X^{1/2-\epsilon}} \right),
\end{align*}
where the third equality is an immediate consequence of Theorem 2.1 of \cite{LOT} in the degree 2 case, for example.  
\end{proof}
\begin{lemma}[Variance of $\omega_Y(n)$]\label{omegavar}
Let $Y \le X^{1/2 - \delta}$ for some $\delta > 0$ and $\mu = \sum_{p \le Y} \frac{1}{p+1}.$
Then
$$\frac{1}{N(X)}\sideset{}{^\flat}\sum_{|D|<X} (\omega_Y(D) - \mu)^2 =\sum_{p \le Y}\frac{1}{p+1}\left(1 - \frac{1}{p+1}\right) + O\left(\frac{Y}{X^{1/2 - \epsilon}}\right).$$
\end{lemma}

\begin{proof}
To calculate the variance, we need to evaluate
\begin{equation} \label{varomegaeq}
\frac{1}{N(X)}\sideset{}{^\flat}\sum_{|D|<X} (\omega_Y(D) - \mu)^2.
\end{equation}

Expanding and distributing the sum over admissible $D$,  this quantity is equal to 
\begin{equation}\label{varexpanded}
\frac{1}{N(X)}\sideset{}{^\flat}\sum_{|D|<X}\left(\sum_{\substack{p < Y \\ p \mid D}} 1\right)^2 - \frac{2 \mu}{N(X)} \sideset{}{^\flat}\sum_{|D|<X}\omega_Y(D) +  \mu^2.
\end{equation}
We'll address the two sums in (\ref{varexpanded}) individually. Note that sums are taken over primes and that the notation is consistent with that in \cite{LOT}.  The leftmost term of (\ref{varexpanded}) give us that
\begin{align}
\frac{1}{N(X)}\sideset{}{^\flat}\sum_{|D|<X}\left(\sum_{\substack{p < Y\\ p \mid D}}1\right)^2 &= \frac{1}{N(X)}\sideset{}{^\flat}\sum_{|D|<X} \left(\sum_{\substack{p,q <Y \\ p,q \mid D}} 1\right)\nonumber\\
&= \frac{1}{N(X)} \sum_{p,q <Y}\sideset{}{^\flat}\sum_{\substack{|D|<X \\ p,q \mid D}}1 \nonumber\\
&= \frac{1}{N(X)}\sum_{p=q < Y} \sideset{}{^\flat}\sum_{\substack{|D|<X\\ p\mid D}} 1 + \frac{1}{N(X)}\sum_{\substack{p,q < Y \\ p\neq q}} \sideset{}{^\flat}\sum_{\substack{|D|<X \\ p,q\mid D}}1\nonumber \\
&= \sum_{p \le Y}\frac{1}{p+1} + \sum_{\substack{p,q <Y \\ p\neq q}}\left( \frac{1}{(p+1)(q+1)} + O(X^{-1/2 + \epsilon}(pq)^{-1/2})\right) \nonumber\\
&= \sum_{p \le Y}\frac{1}{p+1} + \sum_{\substack{p,q <Y \\ p\neq q}}\frac{1}{(p+1)(q+1)} +O\left(\frac{Y}{X^{1/2 - \epsilon}}\right).\label{varpt1}
\end{align}
Where the last line follows from Theorem 2.1 of \cite{LOT}. Now we address the  second term of (\ref{varexpanded}). It gives us that
\begin{align}
\frac{2 \mu}{N(X)}\sideset{}{^\flat}\sum_{|D|<X} \omega_Y(D) &= 2\left(\sum_{p \le Y}\frac{1}{p+1}\right)\left(\sum_{q \le Y}\left[\frac{1}{q+1} + O\left(X^{-1/2 + \epsilon}q^{-1/2}\right)\right] \right) \nonumber  \\
& = 2\left(\sum_{p,q \le Y} \frac{1}{(p+1)(q+1)} + O\left(\frac{Y^{1/2}\log \log Y}{X^{1/2 - \epsilon}}\right) \right). \label{varpt2}
\end{align}
Substituting (\ref{varpt1}) and (\ref{varpt2}) into (\ref{varexpanded}) and expanding $\mu^2$, we find
\begin{equation*}
\frac{1}{N(X)}\sideset{}{^\flat}\sum_{|D|<X} (\omega_Y(D) - \mu)^2 =\sum_{p \le Y}\frac{1}{p+1} - \sum_{p \le Y}\frac{1}{(p+1)^2} + O\left(\frac{Y}{X^{1/2 - \epsilon}} \right).
\end{equation*}
\end{proof}
Now that we have some basic statistical facts about $\omega_Y(D)$, we can use Chebyshev's inequality to give a bound for $\omega_Y(D)$ for most $D$. For our purposes, we will take $Y = \log X$, which will give us mean and variance $\log \log \log X + O(1).$

\begin{theorem}\label{exceptionaldisc}
Let $Y = \log X$. For all $\epsilon >0$, all but $O\left(\dfrac{X}{(\log \log \log X)^{1-\epsilon}}\right)$ fundamental discriminants $D$ with $|D| < X$ are such that $\omega_Y(D) \geq \log \log \log X - O\left((\log \log \log X)^{1 -\epsilon}\right)$.
\end{theorem}

\begin{proof}
Use Chebyshev's inequality with Lemmas \ref{omegamean} and \ref{omegavar} taking $k = (\log \log \log X)^{1/2 - \epsilon}$.
\end{proof}

Applying Theorem \ref{exceptionaldisc} to (\ref{orderofsecond}) gives us that for 100\% of fundamental discriminants $D$, the second moment of the family of quadratic characters modulo $|D|$ is $O\left((\log\log|D|)^{2-\log 2 + \epsilon}\right).$ Using this we can now prove Theorem \ref{typicalchi}. 

\begin{proof}[Proof of Theorem 1.4]
Let $\epsilon, \delta > 0$ and let $D$ be a fundamental discriminant such that $|\omega_Y(D) - \log \log \log |D|| $ is $  O\left((\log\log\log |D|)^{1 - \epsilon}\right).$ By Theorem \ref{exceptionaldisc}, 100\% of fundamental disciriminants have this property.

Lemma \ref{secmoment} above shows that the second moment of the quadratic characters modulo $|D|$ is at most $O\left((\log \log |D|)^{2 - \log 2 + \epsilon} \right),$ when we take $Y = \log D.$ 
Choosing $k = \delta^{-1/2}$ in Corollary \ref{generallowerL}, the result follows.
\end{proof}

\section{Proof of Main Theorem}\label{mainproof}
We are now ready to address ourselves to the proof of Theorem \ref{maintheorem}. The main idea is to use Theorem \ref{typicalchi} to control the residues of the residues of Dedekend $\zeta$-functions appearing in the $D_4$ and $S_4$ estimates in ($\ref{d4s4bounds}$). 
\begin{proof}[Proof of Theorem  \ref{maintheorem}.]

Let $F=\Q(\sqrt{d})$, with $d$ squarefree, be a quadratic number field. Let $W$ be the set quadratic extensions $L/F$ where $L=\Q(\sqrt{d_1}, \sqrt{d_2})$ and $d_1d_2=d$. Note that these extensions are such that $\zeta_L(s)$ factors in such a way as to give us (\ref{residueratio}), and so
\begin{equation}\label{quadraticestimate}
\lim_{X \to \infty} \frac{N_4^F(X, D_4)}{N_4^F(X, S_4)} \gg  \sum_{\substack{L \in W \\ [L:F]=2 }}\frac{\text{Res}_{s=1} \zeta_L(s)}{\text{Res}_{s=1}\zeta_F(s)}.
\end{equation} 
This will be our starting point.  First we'll estimate the residue term by combining (\ref{residueratio}) with Theorem \ref{typicalchi} to see that for $100\%$ quadratic extensions $F$ and sufficiently small $\epsilon$ and $\delta$, and a constant $c$, a proportion $\frac{1-\delta}{2}$ of $L\in W$ satisfy,
$$\frac{\res \zeta_L(s)}{\res \zeta_F(s)} \gg_{\delta} \exp(-c(\log \log |D_F|)^{1-\frac{\log 2}{2} + \epsilon})(\log \log |D_F|)^{-2},$$
the $\log \log |D_F|$ appearing on the right above being a correction factor created when we pass from quadratic characters modulo $|D_F|$ to the $L$-functions of the the corresponding quadratic number fields. 

For any $F$, $\#W = 2^{\omega(d)}$. From Section 2.3 of \cite{Monty} and Theorem \ref{exceptionaldisc} we have that $100\%$ of quadratic fields $F=\Q(\sqrt{d})$ are such that $\#W=2^{\omega(d)} \gg (\log|D_F|)^{\log 2 -\epsilon'}$ for any $\epsilon' > 0$.

So for $100\%$ of quadratic fields $F$, we conclude that
\begin{equation}\label{quadraticestimate2}
\lim_{X \to \infty} \frac{N_4^F(X, D_4)}{N_4^F(X, S_4)} \gg_{\delta}  (\log|D_F|)^{\log 2 -\epsilon'}\exp(-(\log \log |D_F|)^{1-\frac{\log 2}{2} + \epsilon})(\log \log |D_F|)^{-2}.\end{equation}
The statement of Theorem \ref{maintheorem} follows.
\end{proof}

It is worth remarking that we have invoked two \textit{different} $100\%$ results above. The first is that $100\%$ of quadratic number fields are such the bounds on $\omega_Y(D_F)$ are met (in order to get the bound on residues).  The other is that $100\%$ of quadratic number fields are such that the right condition on $\omega(D_F)$ is met (in order to get a suitable bound on $2^{\omega(D_F)}$). In the worst case, the exceptional sets for each of these results are distinct, but their proportion still goes to zero as our bound $X$ on admissible $D_F$ grows. 

Applying the same reasoning from the proof for Theorem \ref{maintheorem}, but without using Theorem \ref{exceptionaldisc}, gives the following lower bound on the ratio for any quadratic number field F.

\begin{theorem}\label{nonstat}
Let $F$ be a quadratic number field and fix $Y \leq |D_F|$, then there is some constant $c$ such that
$$ \lim_{X \to \infty}\frac{N_4^F(X; D_4)}{N_4^F(X; S_4)} \gg \frac{\#\mathrm{Cl}_F[2]}{(\log \log |D_F|)^2(\log |D_F|)^{c\sqrt{ \frac{\log \omega_Y(D_F) + \log \log Y}{\log \log |D_F|} + \frac{1}{2^{\omega_Y(D_F)}}  }}}. $$
\end{theorem}

For a choice of quadratic number field $F$, the above expression shows that if you can pick $Y$ sufficiently small such that $\omega_Y(D_F)$ is sufficiently large, there will be a bias in favor of $D_4$ quartic extensions of $F$. Further, if $\omega_Y(D_F)$ is very large, one should expect that $\#\mathrm{Cl}_F[2]$ is large as well. For example,  if we have $\omega_Y(D_F)$ is of size $\frac{\log |D_F|}{\log \log |D_F|}$ then $\#\textrm{Cl}_F[2]$ is at least of size $|D_F|^{\log 2 / \log\log |D_F|}$.

\section{Examples}\label{examples}

Now that we've shown that most quadratic number fields have more $D_4$ quartic extensions than $S_4$, a couple of natural problems to address are constructing an explicit family of quadratic number fields with arbitrarily more $D_4$ than $S_4$ extensions, and finding the first quadratic number field with more $D_4$ than $S_4$ extensions.

For the first question, consider the family of number fields obtained by taking $F = \Q(\sqrt{\pm d})$ where $d = \prod_{p \le y} p,$ as we take $y \to \infty$. For this family,
$$\omega(D_F) = \frac{\log |D_F|}{\log \log |D_F|}(1 + O(1/\log\log|D_F|)),$$
which immediately gives that $\#\mathrm{Cl}_F[2]$ is about $\exp\left(\frac{\log 2 \log|D_F|}{\log\log |D_F|} \right).$ Because $\omega(D_F)$ is larger than average in this case, we can show that fields in this family have arbitrarily more $D_4$ than $S_4$ extensions without appealing to Theorem \ref{typicalchi}. Instead we can use a lower bound on $L(1,\chi)$ given by Theorem 11.4 in \cite{Monty} which is only conditional on $\chi$ not having an exceptional zero.

In fact, because the formulae from \cite{BSW,CDO} are explicit, we can effectively approximate the constants $\lim_{X \to \infty}\dfrac{N_4^F(X;D_4)}{X}$ and $\lim_{X \to \infty}\dfrac{N_4^F(X;S_4)}{X}$, using either Sage or Magma. In Table \ref{t1}, we see that in our family of fields the percentage of $D_4$ extensions quickly exceeds the percentage of $S_4$ extensions.

Using the same code, we can also answer the question of which quadratic number field is the ``first'' one with more $D_4$ extensions than $S_4$ extensions. Again, we assume $F = \Q(\sqrt{\pm d})$, but now $d$ runs over square-free numbers rather than only the productive of all primes up to $y$ as above. If we order by $|d|$, then we see that about 56\% of quartic extensions of $\Q(\sqrt{-10})$ are $D_4$. See Table \ref{t2}.

\begin{table}
\begin{center}
\begin{tabular}{|l|l|l|l|} 
\hline
$\pm d$     & $S_4$ \textbf{Constant} & $D_4$ \textbf{Constant} & $D_4$ \textbf{Percentage}  \\ 
\hline
2       & 0.06125                 & 0.00255                 & 3.99445                    \\ 
\hline
-2      & 0.02868                 & 0.00242                 & 7.77024                    \\ 
\hline
6       & 0.09898                 & 0.03626                 & 26.81255                   \\ 
\hline
-6      & 0.03389                 & 0.03049                 & 47.35530                   \\ 
\hline
30      & 0.12119                 & 0.20786                 & 63.16992                   \\ 
\hline
-30     & 0.02911                 & 0.11788                 & 80.19609                   \\ 
\hline
210     & 0.11894                 & 0.68112                 & 85.13409                   \\ 
\hline
-210    & 0.02161                 & 0.26399                 & 92.43194                   \\ 
\hline
2310    & 0.13033                 & 1.95228                 & 93.74184                   \\ 
\hline
-2310   & 0.02662                 & 0.75727                 & 96.60405                   \\ 
\hline
30030   & 0.08761                 & 3.14195                 & 97.28722                   \\ 
\hline
-30030  & 0.02961                 & 1.81818                 & 98.39736                   \\ 
\hline
510510  & 0.11305                 & 8.63748                 & 98.70812                   \\ 
\hline
-510510 & 0.02499                 & 3.27599                 & 99.24306                   \\
\hline
\end{tabular}
\end{center}
\caption{}
\label{t1}
\end{table}

\begin{table}
\begin{center}
\begin{tabular}{|l|l|l|l|} 
\hline
\textbf{$\pm d$} & \textbf{$S_4$ Constant} & \textbf{$D_4$ Constant} & \textbf{$D_4$ Percentage}   \\ 
\hline
-1                    & 0.01916                 & 0.00080                 & 4.00075                    \\ 
\hline
2                     & 0.06125                 & 0.00241                 & 3.77973                    \\ 
\hline
-2                    & 0.02868                 & 0.00235                 & 7.55794                    \\ 
\hline
3                     & 0.07729                 & 0.02138                 & 21.66628                   \\ 
\hline
-3                    & 0.01480                 & 0.00015                 & 1.01581                    \\ 
\hline
5                     & 0.04181                 & 0.00041                 & 0.97732                    \\ 
\hline
-5                    & 0.03783                 & 0.02618                 & 40.90038                   \\ 
\hline
6                     & 0.09898                 & 0.03602                 & 26.68166                   \\ 
\hline
-6                    & 0.03389                 & 0.03025                 & 47.16238                   \\ 
\hline
7                     & 0.11253                 & 0.03552                 & 23.99301                   \\ 
\hline
-7                    & 0.02954                 & 0.00051                 & 1.68833                    \\ 
\hline
10                    & 0.12577                 & 0.07665                 & 37.86747                   \\ 
\hline
-10                   & 0.02468                 & 0.03141                 & 55.99729  \\
\hline
\end{tabular}
\end{center}
\caption{}
\label{t2}
\end{table}

\end{document}